\documentclass[11pt]{article}

\usepackage[a4paper,left=2cm,right=2cm,top=1.5cm,bottom=2cm]{geometry}
\usepackage[utf8]{inputenc}
\usepackage[T1]{fontenc}
\usepackage[english]{babel}
\usepackage{amsthm}
\usepackage{amsmath}
\usepackage{amssymb}
\usepackage{mathrsfs}
\usepackage{setspace}
\usepackage{float}
\usepackage{multirow}
\usepackage{stmaryrd}
\usepackage{color}
\usepackage{wrapfig}
\usepackage{pifont}
\usepackage{array}
\usepackage[colorlinks=true,linkcolor=ocre,citecolor=ocre]{hyperref}
\usepackage{dsfont}
\usepackage{upgreek}
\usepackage{nicefrac}
\usepackage{tikz}{}
\usepackage{gensymb}
\usepackage{FiraSans}
\usepackage{graphicx}
\usepackage{caption}
\renewenvironment{proof}{\noindent{\sffamily{\textbf{Proof :}}}}{\begin{flushright}$\square$\end{flushright}}

\newcommand{\IE}{\mathbb{E}}
\newcommand{\IN}{\mathbb{N}}
\newcommand{\IZ}{\mathbb{Z}}

\newcommand{\IR}{\mathbb{R}}

\newcommand{\IT}{\mathbb{T}}

\newcommand{\drm}{\mathrm d}

\newcommand{\CD}{\mathcal D}
\newcommand{\CS}{\mathcal S}

\newcommand{\CC}{\mathcal C}
\newcommand{\CH}{\mathcal H}

\newcommand{\CX}{\mathcal X}

\newcommand{\CB}{\mathcal B}

\newcommand{\SF}{\mathscr{F}}

\newcommand{\IDC}{\mathds{1}}

\newcommand{\eps}{\varepsilon}

\definecolor{ocre}{RGB}{64,123,121}

\newcounter{item}
\numberwithin{item}{section}

\newtheorem{theorem}[item]{\sffamily Theorem}
\newtheorem{definition}[item]{\sffamily Definition}
\newtheorem{proposition}[item]{\sffamily Proposition}

\newtheorem*{theorem*}{\sffamily Theorem}
\newtheorem*{definition*}{\sffamily Definition}
\newtheorem*{proposition*}{\sffamily Proposition}
\newtheorem*{lemma*}{\sffamily Lemma}
\newtheorem*{corollary*}{\sffamily Corollary}

\providecommand{\keywords}[1]
{
	{\footnotesize	
	\textbf{Keywords --} #1}
}

\usepackage[explicit]{titlesec}
\titleformat{\section}{\centering\Large\bfseries}{\thesection \ --}{0.7em}{\Large\bfseries #1}
\titleformat{\subsection}{\centering\large\bfseries}{\thesubsection \ --}{0.4em}{\large\bfseries #1}
\titleformat{\subsubsection}{\centering\bfseries}{\thesubsubsection \ --}{0.4em}{\bfseries #1}

\let\emph\relax
\DeclareTextFontCommand{\emph}{\bfseries\em}

\title{\bfseries A simple construction of the Anderson operator via its quadratic form in dimensions two and three}
\author{Antoine MOUZARD and El Maati OUHABAZ}
\date{}

\begin{document}

\maketitle
\abstract{We provide a simple construction of the Anderson operator in dimensions two and three. This is done through its quadratic form. We rely on an exponential transform instead of the regularity structures or paracontrolled calculus which are usually used for  the construction of the  operator. The knowledge of the form is robust enough to deduce important properties such as positivity and irreducibility of the corresponding semigroup. The latter property gives  existence of a spectral gap.}
\vspace{0.5cm}


\section*{Introduction}

Over the last decade, the study of singular stochastic PDEs has grown to an important field with the introduction of regularity structures by Hairer \cite{Hai14} and paracontrolled calculus by Gubinelli, Imkeller and Perkowski \cite{GIP}. The theory first aimed at the resolution of parabolic equations such as the Parabolic Anderson Model (PAM) equation or the Kardar–Parisi–Zhang (KPZ) equation, it then led to the construction of the Anderson Hamiltonian 
\begin{equation*}
H=-\Delta+\xi
\end{equation*}
with $\xi$ the spatial white noise, see \cite{AllezChouk,GUZ,Labbe,CVZ,Mouzard} in dimension $2$ and $3$, on a finite box with periodic or Dirichlet boundary conditions or even compact Riemaniann manifolds. 

In this note we provide a simple construction of this operator via its quadratic form without using regularity structures or paracontrolled calculus. We rely on an exponential transform first used by Hairer and Labbé \cite{HairerLabbe15} for the continuum parabolic Anderson model on $\IR^2$ and then used in different context, see for example \cite{DW,TV,JagannathPerkowski23,BDFT23one,ChauleurMouzard23}. In particular, this was already used by Matsuda and van Zuijlen \cite{MatsudaZuijlen22} to construct the Anderson form in the full subcritical regime using also regularity structures. See also \cite{MM,Mouzard22} for other singular stochastic operators.

The Anderson Hamiltonian is the Schrödinger operator $H=-\Delta+\xi$ with $\xi$ the space white noise which is a random distribution of negative Hölder regularity $-\frac{d}{2}-\kappa$ for any $\kappa>0$. In one dimension, it is the derivative of the Brownian motion and the associated form
\begin{equation*}
\mathfrak{a}(u,v)=\int_0^1\nabla u(x)\cdot\nabla v(x)\drm x+\int_0^1u(x)v(x)\xi(\drm x)
\end{equation*}
was constructed by Fukushima and Nakao \cite{FN} with domain the usual Sobolev space $\CH^1$. The idea is that one can multiply two distributions if the sum of their regularity is positive, hence $uv\xi$ is well-defined as a distribution for $u,v\in\CH^1$ since $\xi\in\CC^{-\frac{1}{2}-\kappa}$. In two dimensions, $\xi\in\CC^{-1-\kappa}$ and this construction is not possible anymore. Following the recent progress in singular stochastic PDEs, the operator can be constructed with a random domain $\CD_\xi^2$ depending on the noise $\xi$ such that $H:\CD_\xi^2\subset L^2\rightarrow L^2$ is an unbounded closed operator. Taking $u\in L^2$ and assuming that $Hu$ is an element of $L^2$, one obtains  the relation
\begin{equation*}
\Delta u=u\xi-Hu
\end{equation*}
which induces an expansion of $u$ with respect to the noise using regularity structures or paracontrolled calculus. In particular, \cite{GUZ} and \cite{Mouzard} also identify a form domain, that is a random subspace $\CD_\xi^1\subset L^2$ such that
\begin{equation*}
\forall u\in\CD_\xi^1,\quad|\langle Hu,u\rangle|<\infty.
\end{equation*}
We emphasize in the notation that the domains are random and consist of random functions depending on the noise. In each case, the operator is a singular stochastic operator and a renormalization procedure is involved in its construction. For a regularization of the noise $\xi_\eps$, the operator is constructed as a limit in the resolvent sense, that is
\begin{equation*}
H=\lim_{\eps\to0}\big(-\Delta+\xi_\eps-c_\eps\big)
\end{equation*}
with a constant $c_\eps$ which explodes when $\eps \to 0$. It is related to the definition of the product $\Delta^{-1}\xi\cdot\xi$ and the divergence of the Green function of the Laplacian.  In two dimensions one has $c_\eps\sim\log(\eps)$ while $c_\eps\sim\eps^{-1}$ in three dimensions. In one dimension, this product is well-defined and one can take $c_\eps=0$ which is coherent with \cite{FN} and $\CD_\xi^1=\CH^1$ does not depend on the noise. However the domain of the operator is random and this method was recently used by Dumaz and Labbé to provide a precise study of the operator, see \cite{DumazLabbe22} and references therein. The Anderson form also appears as the energy for dispersive PDEs such as the nonlinear Schrödinger equation
\begin{equation*}
i\partial_tu=\Delta u+u\xi+|u|^2u
\end{equation*}
and was used to obtain solutions, see for example \cite{DW,GUZ,Mouzard,TV,MZ,ChauleurMouzard23} and references therein. In particular, uniform bounds in energy is the crucial property of such singular stochastic PDEs where one does not have the regularizing properties of the parabolic equation. In this context, one has to work with random initial data depending on the noise and the conservation of energy makes the form domain of the Anderson operator a natural space to get a global solution.

In this work, we consider a new variable $u=e^Xv$ for a suitable random field $X$. In this case, we have
\begin{equation*}
\Delta u=e^X\Delta v+2e^X\nabla X\cdot\nabla v+e^X(|\nabla X|^2+\Delta X)v
\end{equation*}
and if $X$ is a solution to $\Delta X=\xi$, the Anderson operator is formaly given by
\begin{equation*}
Hu=-e^X\Delta v-2e^X\nabla X\cdot\nabla v-e^X|\nabla X|^2v.
\end{equation*}
In two dimensions, we have
\begin{equation*}
\xi\in\CC^{-1-\kappa}\quad\implies\quad\nabla X\in\CC^{-\kappa},
\end{equation*}
hence the square $|\nabla X|^2$ is singular and has to be defined with a renormalization procedure as a Wick product $|\nabla X|^{2\diamond}\in\CC^{-2\kappa}$. In this case, $v\in\CH^1$ is regular enough for the associated form to make sense and one can construct the Anderson form with domain $\CD(\mathfrak{a})=e^X\CH^1$. In three dimensions, we have
\begin{equation*}
\xi\in\CC^{-\frac{3}{2}-\kappa}\quad\implies\quad\nabla X\in\CC^{-\frac{1}{2}-\kappa}
\end{equation*} 
and the Wick product $|\nabla X|^{2\diamond}\in\CC^{-1-2\kappa}$ is too rough to be multiplied by $v\in\CH^1$. One can apply the same method and construct the Anderson form with domain $e^{X+Y}\CH^1$ with a suitable second random field $Y$.

This exponential transform  allows us to construct a symmetric form $\mathfrak{a}$ whose associated  operator $H$  is the Anderson Hamiltonian. The domain of $H$ is given by
\begin{equation*}
\CD(H)=\big\{u\in L^2(\IT^d)\ ;\ \exists v\in L^2(\IT^d),\forall\varphi\in\CD(\mathfrak{a}), \mathfrak{a}(u,\varphi)=\langle v,\varphi\rangle\big\}.
\end{equation*}
For a better description of this domain a more involved theory such as regularity structures or paracontrolled calculus seems to be  necessary. Nevertheless, the knowledge  of $H$ through its form is enough to deduce that $H$ is self-adjoint, it has a discrete spectrum and an $L^2$ orthonormal basis given by eigenfunctions. In addition, relying on a criterion from  \cite{Ouhabaz05},  we prove that the associated semigroup is irreducible. In particular, this implies the existence of a spectral gap $\lambda_1<\lambda_2$ with a positive ground state $\Psi\in\CD(H)$. This result was already proved in \cite{BailleulDangMouzard22} by relying on a quantitative estimate for the linear Parabolic Anderson Model equation. Our work provides a pedestrian approach to this result even in three dimensions which usually relies on involved computations with expansion of order $5$ using regularity structures or paracontrolled calculus.

In order to keep the ideas and tools simple we restrict ourselves  to the case of the torus $\IT^d$ for $d \in \{2,3\}$ (endowed with the Lebesgue measure $dx$).  Our construction works on any compact manifold without boundary. Following the tools used in \cite{Mouzard} we could also allow  the manifold to have a smooth boundary and the operator is then subject to  Dirichlet boundary conditions. The approach using the exponential transform  can not be used to construct the explicit domain of the operator, the use of regularity structure or paracontrolled calculus seems inevitable. 

In Section \ref{SectionStoc}, we provide several    bounds on stochastic functions and distributions that we need to construct the form. In particular, this is where the renormalization of the singular products is done. In Sections \ref{Section2D} and \ref{Section3D}, we respectively construct the Anderson form in two and three dimensions using the first and second order exponential transform. In Section \ref{SectionSpectralGap}, we prove irreducibility and existence of a spectral gap. 

\section{Stochastic bounds and renormalization} \label{SectionStoc}

On the torus $\IT^d$, the white noise is given formally by
\begin{equation*}
\xi(x)=\sum_{k\in\IZ^d}\xi_ke^{ik\cdot x}
\end{equation*}
with $(\xi_k)_{k\in\IZ^d}$ a family of independent and identically distributed random variables of centered standard complex Gaussian with $\xi_{-k}=\overline{\xi_k}$. This gives a centered real Gaussian field with covariance function
\begin{equation*}
\IE\big[\xi(x)\xi(y)\big]=\delta_0(x-y),
\end{equation*}
that is a random distribution $(\langle\xi,\varphi\rangle)_{\varphi\in L^2(\IT^d)}$ such that
\begin{equation*}
\IE\big[\langle\xi,\varphi\rangle\langle\xi,\psi\rangle\big]=\langle\varphi,\psi\rangle_{L^2(\IT^d)}
\end{equation*}
which indeed gives the Fourier coefficient $(\xi_k)_{k\in\IZ^d}$. Its first construction is due to Paley and Zygmund \cite{PaleyZygmund30,PaleyZygmund32} and is actually the first random distribution ever considered. A natural and convenient setting is given by the Besov space $\CB_{p,q}^\alpha$ which can be defined using the Littlewood-Paley  decomposition, see for example \cite{BCD}. This decomposition can be stated  as follows 
\begin{equation*}
u=\sum_{n\ge0}\Delta_nu
\end{equation*}
with $\Delta_nu=\big(\SF^{-1}\IDC_{|\cdot|\simeq 2^n}\SF\big)u$, that is the projection of $u$ in frequencies on an annulus of size $2^n$. It is defined by 
\begin{equation*}
\big(\Delta_nu\big)(x):=2^{d(n-1)}\int_{\IR^d}\chi\big(2^{n-1}(x-y)\big)u(y)\drm y
\end{equation*} 
with $\chi\in\CS(\IR^d)$ and $\text{supp}\ \widehat\chi\subset\{\frac{1}{2}\le|z|\le 2\}$ for $n\ge1$ and
\begin{equation*}
\big(\Delta_0u\big)(x):=\int_{\IR^d}\chi_0(x-y)u(y)\drm y
\end{equation*}
with $\chi_0\in\CS(\IR^d)$ and $\text{supp}\ \widehat\chi_0\subset\{|z|\le 1\}$. We also denote $K=\widehat{\chi}$ such that $\Delta_nu=\big(\SF^{-1}K(2^n\cdot)\SF\big)u$. Then the Besov space $\CB_{p,q}^\alpha$ are distributions such that
\begin{equation*}
\|u\|_{\CB_{p,q}^\alpha}:=\Big(\ \sum_{n\ge0}2^{\alpha pn}\|\Delta_nu\|_{L^q(\IT^d)}^p\ \Big)^{\frac{1}{p}}<\infty.
\end{equation*}
The particular case $p=q=2$ corresponds to the Sobolev space $\CB_{2,2}^\alpha=\CH^\alpha$ and for $p=q=\infty$ with $\alpha\in\IR_+\backslash\IN$, one gets the usual Hölder spaces $\CB_{\infty,\infty}^\alpha=\CC^\alpha$. One also has the continuous Besov embedding
\begin{equation*}
\CB_{p_1,q_1}^\alpha\hookrightarrow \CB_{p_2,q_2}^{\alpha-d(\frac{1}{p_1}-\frac{1}{p_2})}
\end{equation*}
for $p_1\le p_2,q_1\le q_2$ and $\alpha\in\IR$. While one can a priori only multiply a distribution by a smooth function, one has the following product rule in the case of Besov spaces which corresponds to Young condition.

\medskip

\begin{proposition}\label{PropProduct}
For $\alpha,\beta\in\IR$ such that $\alpha+\beta>0$ and $p,q,r\in[1,\infty]$ such that $\frac{1}{r}=\frac{1}{p}+\frac{1}{q}$, there exists a constant $C>0$ such that
\begin{equation*}
\|uv\|_{\CB_{r,r}^{\alpha\wedge\beta}}\le C\|u\|_{\CB_{p,p}^\alpha}\|v\|_{\CB_{q,q}^\beta}.
\end{equation*}
\end{proposition}

\medskip

The following proposition gives a similar result at the level of the duality bracket.

\medskip

\begin{proposition}\label{PropDuality}
For $\alpha\in\IR$ and $p,p',q,q'\in[1,\infty]$ such that $1=\frac{1}{p}+\frac{1}{p'}=\frac{1}{q}+\frac{1}{q'}$, there exists a constant $C>0$ such that
\begin{equation*}
|\langle u,v\rangle|\le C\|u\|_{\CB_{p,q}^\alpha}\|v\|_{\CB_{p',q'}^{-\alpha}}.
\end{equation*}
\end{proposition}

\medskip

For later use, we introduce a new random field $X$ defined by 
\begin{equation*}
X(x)=-\sum_{k\in\IZ^d\backslash\{0\}}\frac{1}{|k|^2}\xi_ke^{ik\cdot x}. 
\end{equation*}
It  satisfies the equation
\begin{equation*}
\Delta X=\xi-\xi_0.
\end{equation*}
The following proposition gives  Hölder regularity of $\xi$ and $X$.

\medskip

\begin{proposition}
For any $\kappa>0$, one has almost surely
\begin{equation*}
\xi\in\CC^{-\frac{d}{2}-\kappa}(\IT^d)\quad\text{and}\quad X\in\CC^{2-\frac{d}{2}-\kappa}(\IT^d).
\end{equation*}
\end{proposition}

\medskip

\begin{proof}
Since the noise is Gaussian, we have
\begin{equation*}
\IE\big[\langle \xi,\varphi\rangle^p\big]\le(p-1)^{\frac{p}{2}}\IE\big[\langle\xi,\varphi\rangle^2\big]^{\frac{p}{2}}
\end{equation*}
for any test function $\varphi$.  This is usually referred to as Gaussian hypercontractivity. In order to use this, we estimate the Besov norm $\CB_{p,p}^\gamma$ for $p$ large and use the embedding
\begin{equation*}
\CB_{p,p}^\gamma(\IT^d) \hookrightarrow  \CB^{\gamma-\frac{d}{p}}_{\infty,\infty}(\IT^d).
\end{equation*}
We have
\begin{align*}
\IE\Big[\|\Delta_n\xi\|_{L^p(\IT^d)}^p\Big]&=\int_{\IT^d}\IE\big[\langle \xi,\chi_n(x-\cdot)\rangle^p\big]\drm x\\
&\le(p-1)^{\frac{p}{2}}\int_{\IT^d}\IE\big[\langle \xi,\chi_n(x-\cdot)\rangle^2\big]^{\frac{p}{2}}\drm x\\
&\le(p-1)^{\frac{p}{2}}\|\chi_n\|_{L^2(\IT^d)}^p|\IT^d|
\end{align*}
with $\chi_n(\cdot)=2^{dn}\chi(2^n\cdot)$ and using that $\xi$ is an isometry from $L^2(\IT^d)$ to $L^2(\Omega)$.  We have
\begin{equation*}
\|\chi_n\|_{L^2(\IT^d)}^2=2^{2dn}\|\chi(2^n\cdot)\|_{L^2(\IT^d)}^2=2^{dn}\|\chi\|_{L^2(\IT^d)}^2
\end{equation*}
hence
\begin{equation*}
\IE\Big[\|\Delta_n\xi\|_{L^p(\IT^d)}^p\Big]\le(p-1)^{\frac{p}{2}}2^{pn\frac{d}{2}}\|\chi\|_{L^2(\IT^d)}^p|\IT^d|.
\end{equation*}
This gives
\begin{equation*}
\IE\big[\|\xi\|_{\CB_{p,p}^{-\frac{d}{2}}}\big]<\infty,
\end{equation*}
and hence $\xi\in\CC^{-\frac{d}{2}-\frac{d}{p}}(\IT^d)$ for any $p\ge1$ which completes the proof for the regularity of $\xi$ while the regularity of $X$ follows from a standard regularity estimate.
\end{proof}

In two dimensions, one has $X\in\CC^{1-\kappa}$ hence $\nabla X\in\CC^{-\kappa}$ and the square $|\nabla X|^2$ is ill-defined since $-2\kappa<0$. Consider a regularization of the noise $\xi_\eps=\xi*\rho_\eps$ with $\rho_\eps$ a mollifier. Then $\xi_\eps$ converges to $\xi$ as $\eps$ goes to $0$ in $\CC^{-1-\kappa}$ and one can consider $X_\eps$ the solution to
\begin{equation*}
\Delta X_\eps=\xi_\eps-\langle\xi_\eps,1\rangle
\end{equation*}
which converges to $X$ in $\CC^{1-\kappa}$ as $\eps$ goes to $0$. Since the square $|\nabla X|^2$ is ill-defined, the quantity $|\nabla X_\eps|^2$ diverges and this is described by the Wick square as proved in the following proposition.

\medskip

\begin{proposition}
There exists a distribution $|\nabla X|^{2\diamond}\in\CC^{-2\kappa}(\IT^2)$ such that
\begin{equation*}
|\nabla X|^{2\diamond}=\lim_{\eps\to0}\Big(|\nabla X_\eps|^2-\IE\big[|\nabla X_\eps|^2\big]\Big)
\end{equation*}
in $\CC^{-2\kappa}(\IT^2)$ in probability. Moreover, one has
\begin{equation*}
\IE\big[|\nabla X_\eps|^2\big]\underset{\eps\to0}{\sim}\frac{1}{(2\pi)^2}\log(\eps).
\end{equation*}
\end{proposition}

\medskip

\begin{proof}
Since $\xi_\eps=\xi*\rho_\eps$, we have
\begin{equation*}
X_\eps(x)=-\sum_{k\in\IZ^2\backslash\{0\}}\frac{\widehat{\rho_\eps}(k)}{|k|^2}\xi_ke^{ikx}
\end{equation*}
thus
\begin{equation*}
|\nabla X_\eps(x)|^2=\sum_{k,k'\in\IZ^2\backslash\{0\}}\widehat{\rho_\eps}(k)\overline{\widehat{\rho_\eps}(k')}\frac{k\cdot k'}{|k|^2|k'|^2}\xi_k\overline{\xi_{k'}}e^{i(k-k')\cdot x}.
\end{equation*}
Using $\IE\big[\xi_k\overline{\xi_{k'}}\ \big]=\delta_0(k-k')$, we have 
\begin{equation*}
\IE\big[|\nabla X_\eps(x)|^2\big]=-\sum_{k\in\IZ^2\backslash\{0\}}\frac{|\widehat{\rho_\eps}(k)|^2}{|k|}\underset{\eps\to0}{\sim}\frac{1}{(2\pi)^2}\log(\eps)
\end{equation*}
which gives the second part of the statement. For $n\ge1$, we have
\begin{equation*}
\Delta_n\Big(|\nabla X_\eps(x)|^2-\IE\big[|\nabla X_\eps(x)|^2\big]\Big)=\sum_{k\neq k'}K_n(k-k')\widehat{\rho_\eps}(k)\overline{\widehat{\rho_\eps}(k')}\frac{k\cdot k'}{|k|^2|k'|^2}\xi_k\overline{\xi_{k'}}e^{i(k-k')\cdot x}
\end{equation*}
with $K_n(\cdot)=K(2^{-n}\cdot)$ and so  $\Big|\Delta_n\Big(|\nabla X_\eps(x)|^2-\IE\big[|\nabla X_\eps(x)|^2\big]\Big)\Big|^2$ is given by
\begin{eqnarray*}
\sum_{\substack{k_1\neq k_1'\\k_2\neq k_2'}}K_n(k_1-k_1')\overline{K_n(k_2-k_2')}\widehat{\rho_\eps}(k_1)\overline{\widehat{\rho_\eps}(k_1')}\overline{\widehat{\rho_\eps}(k_2)}\widehat{\rho_\eps}(k_2')\frac{k_1\cdot k_1'}{|k_1|^2|k_1'|^2}\frac{k_2\cdot k_2'}{|k_2|^2|k_2'|^2}\xi_{k_1}\overline{\xi_{k_1'}}\xi_{k_2}\overline{\xi_{k_2'}}\\
\hspace{1.5cm} \times e^{i(k_1-k_1')\cdot x}e^{-i(k_2-k_2')\cdot x}
\end{eqnarray*}
and we have
\begin{align*}
\IE\big[\xi_{k_1}\overline{\xi_{k_1'}}\xi_{k_2}\overline{\xi_{k_2'}}\big]&=\IE\big[\xi_{k_1}\overline{\xi_{k_1'}}\big]\IE\big[\xi_{k_2}\overline{\xi_{k_2'}}\big]+\IE\big[\xi_{k_1}\xi_{k_2}\big]\IE\big[\overline{\xi_{k_1'}\xi_{k_2'}}\big]+\IE\big[\xi_{k_1}\overline{\xi_{k_2'}}\big]\IE\big[\overline{\xi_{k_1'}}\xi_{k_2}\big]\\
&=\delta_0(k_1-k_1')\delta_0(k_2-k_2')+\delta_0(k_1+k_2)\delta_0(k_1'+k_2')+\delta_0(k_1-k_2')\delta_0(k_1'-k_2)
\end{align*}
for any $k_1,k_1',k_2,k_2'\in\IZ^2$. The term $k_1=k_1'$ and $k_2=k_2'$ corresponds to $\IE\big[|\nabla X_\eps(x)|^2\big]$ hence the restriction of the sum. It follows that 
\begin{align*}
\IE\Big|\Delta_n\Big(|\nabla X_\eps(x)|^2-\IE\big[|\nabla X_\eps(x)|^2\big]\Big)\Big|^2&=2\sum_{k_1,k_2}|K_n(k_1-k_2)|^2\frac{|\widehat{\rho_\eps}(k_1)|^2|\widehat{\rho_\eps}(k_2)|^2}{|k_1|^2|k_2|^2}\\
&=2\sum_{k_1,k_2}\big|K\big(2^{-n}(k_1-k_2)\big)\big|^2\frac{|\widehat{\rho_\eps}(k_1)|^2|\widehat{\rho_\eps}(k_2)|^2}{|k_1|^2|k_2|^2}\\
&=2\sum_k|K(2^{-n}k)|^2\sum_{k_1-k_2=k}\frac{|\widehat{\rho_\eps}(k_1)|^2|\widehat{\rho_\eps}(k_2)|^2}{|k_1|^2|k_2|^2}\\
&=2\sum_k|K(2^{-n}k)|^2\sum_{k_2}\frac{|\widehat{\rho_\eps}(k+k_2)|^2|\widehat{\rho_\eps}(k_2)|^2}{|k+k_2|^2|k_2|^2}\\
&\le C2^{2n}2^{-(2-2\kappa)n}\sum_{k_2}\frac{|\rho_\eps(k_2)|^2}{|k_2|^{2+2\kappa}}
\end{align*}
for any $\kappa>0$ and a constant $C>0$ using the support of $K$. The Gaussian hypercontractivity yields 
\begin{align*}
\IE\Big|\Delta_n\Big(|\nabla X_\eps(x)|^2-\IE\big[|\nabla X_\eps(x)|^2\big]\Big)\Big|^p&\le(p-1)^p\Big(\IE\Big|\Delta_n\Big(|\nabla X_\eps(x)|^2-\IE\big[|\nabla X_\eps(x)|^2\big]\Big)\Big|^2\Big)^{\frac{p}{2}}\\
&\le C2^{\kappa np}. 
\end{align*}
Thus,  $|\nabla X_\eps(x)|^2-\IE\big[|\nabla X_\eps(x)|^2$ is bounded in $\CB_{p,p}^{-\kappa}$ for any $\kappa>0$ and $p\ge1$. Using the embedding  $\CB_{p,p}^{-\kappa}\hookrightarrow\CC^{-\kappa-\frac{d}{p}}$ and a similar bound, one  proves that $\big(|\nabla X_\eps(x)|^2-\IE\big[|\nabla X_\eps(x)|^2\big)_{\eps>0}$ is a Cauchy family in $\CC^{-\kappa}$ for any $\kappa>0$ which completes the proof.
\end{proof}

We define the two dimensional enhanced noise
\begin{equation*}
\Xi=\big(\xi,|\nabla X|^{2\diamond}\big)
\end{equation*}
which belongs to
\begin{equation*}
\CX^\kappa(\IT^2)=\CC^{-1-\kappa}(\IT^2)\times\CC^{-2\kappa}(\IT^2)
\end{equation*}
for any $\kappa>0$. We also have that
\begin{equation*}
\Xi_\eps=\big(\xi_\eps,|\nabla X_\eps|^2-(2\pi)^{-2}\log(\eps)\big)
\end{equation*}
converges to $\Xi$ in $\CX^\kappa(\IT^2)$ for any $\kappa>0$. In three dimensions, one has $X\in\CC^{\frac{1}{2}-\kappa}$ hence this term is even more singular with $-1-2\kappa<0$. The analog of the previous renormalization is the following proposition with a larger divergence. Its proof follows the same path as the previous one.

\medskip

\begin{proposition}
There exist a distribution $|\nabla X|^{2\diamond}\in\CC^{-1-2\kappa}(\IT^3)$ such that
\begin{equation*}
|\nabla X|^{2\diamond}=\lim_{\eps\to0}\Big(|\nabla X_\eps|^2-\IE\big[|\nabla X_\eps|^2\big]\Big)
\end{equation*}
in $\CC^{-1-2\kappa}(\IT^3)$ in probability. Moreover, one has
\begin{equation*}
\IE\big[|\nabla X_\eps|^2\big]\underset{\eps\to0}{\sim}-\frac{1}{(2\pi)^2}\frac{1}{\eps}.
\end{equation*}
\end{proposition}

\medskip

Since the noise is more irregular, $|\nabla X|^{2\diamond}$ is too rough to make sense of its bracket with $\CH^1$ hence we will also need $Y$ the solution to
\begin{equation*}
\Delta Y=|\nabla X|^{2\diamond}-\langle|\nabla X|^{2\diamond},1\rangle
\end{equation*}
which belongs to $\CC^{-2\kappa}(\IT^3)$. Its square is also singular and can be defined as a Wick product, as well as the product $\nabla X\cdot\nabla Y$.

\medskip

\begin{proposition}
There exists a distribution $|\nabla Y|^{2\diamond}\in\CC^{-4\kappa}(\IT^3)$ such that
\begin{equation*}
|\nabla Y|^{2\diamond}=\lim_{\eps\to0}\Big(|\nabla Y_\eps|^2-\IE\big[|\nabla Y_\eps|^2\big]\Big)
\end{equation*}
in $\CC^{-4\kappa}(\IT^3)$ in probability. Moreover, one has
\begin{equation*}
\IE\big[|\nabla Y_\eps|^2\big]\underset{\eps\to0}{\sim}\frac{1}{(2\pi)^2}\log(\eps).
\end{equation*}
There also exists a distribution $\nabla X\diamond\nabla Y\in\CC^{-\frac{1}{2}-3\kappa}(\IT^3)$ such that
\begin{equation*}
\nabla X\diamond\nabla Y=\lim_{\eps\to0}\Big(\nabla X_\eps\cdot\nabla Y_\eps\Big).
\end{equation*}
\end{proposition}

\medskip

We define the three dimensional enhanced noise
\begin{equation*}
\Xi=\big(\xi,|\nabla X|^{2\diamond},|\nabla Y|^{2\diamond},\nabla X\diamond\nabla Y\big)
\end{equation*}
which belongs to
\begin{equation*}
\CX^\kappa(\IT^3)=\CC^{-\frac{3}{2}-\kappa}(\IT^3)\times\CC^{-1-2\kappa}(\IT^3)\times\CC^{-4\kappa}(\IT^3)\times\CC^{-\frac{1}{2}-3\kappa}
\end{equation*}
for any $\kappa>0$. We also have that
\begin{equation*}
\Xi_\eps=\big(\xi_\eps,|\nabla X_\eps|^2-(2\pi)^{-2}\eps^{-1},|\nabla Y_\eps|^2-(2\pi)^{-2}\log(\eps),\nabla X_\eps\cdot\nabla Y_\eps\big)
\end{equation*}
converges to $\Xi$ in $\CX^\kappa(\IT^3)$ for any $\kappa>0$.

\section{Construction in two dimensions} \label{Section2D}

It is tempting to define the   form of the Anderson operator by
\begin{equation*}
{\mathfrak a}(u_1,u_2)=\int_{\IT^2}\nabla u_1(x)\cdot\nabla u_2(x)\drm x+\int_{\IT^2}u_1(x)u_2(x)\xi(\drm x)
\end{equation*}
for any $u_1,u_2\in C^\infty(\IT^2)$. However, this is not a natural object since this form is not closable as shown by the recent progress on singular stochastic operators, which can be guessed from the fact that for $u\in\CH^1$ the form domain of $\Delta$, the product $u\xi$ is ill-defined. For $\xi_\eps=\xi*\rho_\eps$ a regularization of the noise, consider the regularized form
\begin{equation*}
{\mathfrak a}_\eps(u_1,u_2)=\int_{\IT^2}\nabla u_1(x)\cdot\nabla u_2(x)\drm x+\int_{\IT^2}u_1(x)u_2(x)\big(\xi_\eps(x)-c_\eps\big)\drm x
\end{equation*}
with $c_\eps$ the logarithmic diverging constant defined in the previous section. For any fixed $\eps>0$, ${\mathfrak a}_\eps$ is a closed symmetric form with domain $\CH^1$ and we construct a form ${\mathfrak a}$ such that ${\mathfrak a}_\eps$ converges to ${\mathfrak a}$ as $\eps$ goes to $0$. With $X$ the random field constructed in the previous section, we consider the new variable $u=e^Xv$ and define
\begin{equation*}
Hu:=-e^X\Delta v-2e^X\nabla X\cdot\nabla v-e^X|\nabla X|^{2\diamond}v+\xi_0e^Xv
\end{equation*}
for $v\in\CC^\infty$. Since $X\in\CC^{1-\kappa}$ and $|\nabla X|^{2\diamond}\in\CC^{-2\kappa}$, $Hu$ is well-defined as a distribution. The associated form is given by 
\begin{align*}
{\mathfrak a}(u_1,u_2)&=\langle Hu_1,u_2\rangle\\
&=\langle He^Xv_1,e^Xv_2\rangle\\
&=-\langle\Delta v_1,v_2e^{2X}\rangle-2\langle \nabla X\cdot\nabla v_1,v_2e^{2X}\rangle-\langle|\nabla X|^{2\diamond}v_1,v_2e^{2X}\rangle+\xi_0\langle v_1,v_2e^{2X}\rangle\\
&=\int_{\IT^2}\nabla v_1(x)\cdot\nabla v_2(x)e^{2X(x)}\drm x-\langle|\nabla X|^{2\diamond}v_1,v_2e^{2X}\rangle+\xi_0\int_{\IT^2}v_1(x)v_2(x)e^{2X(x)}\drm x
\end{align*}
which is well-defined for $v_1,v_2\in\CH^1$ since
\begin{align*}
\big|\langle|\nabla X|^{2\diamond}e^{2X},v_1v_2\rangle\big|&\le\||\nabla X|^{2\diamond}e^{2X}\|_{\CC^{-\kappa}}\|v_1v_2\|_{\CB_{1,1}^\kappa}\\
&\le\||\nabla X|^{2\diamond}\|_{\CC^{-\kappa}}\|e^{2X}\|_{\CC^{2\kappa}}\|v_1\|_{\CH^{2\kappa}}\|v_2\|_{\CH^{2\kappa}}\\
&\le\||\nabla X|^{2\diamond}\|_{\CC^{-\kappa}}\|e^{2X}\|_{\CC^{1-\kappa}}\|v_1\|_{\CH^1}\|v_2\|_{\CH^1}
\end{align*} 
for $\kappa>0$ small enough using Proposition \ref{PropDuality} and Proposition \ref{PropProduct}.

\medskip

\begin{definition}
The Anderson form is defined by 
\begin{equation*}
{\mathfrak a}(u_1,u_2):=\langle \nabla v_1,\nabla v_2\rangle_{L^2(\IT^2,e^{2X}\drm x)}-\big\langle|\nabla X|^{2\diamond},v_1v_2e^{2X}\big\rangle+\xi_0\langle v_1,v_2\rangle_{L^2(\IT^2,e^{2X}\drm x)}
\end{equation*}
where $v_i=e^{-X}u_i$ with domain $\CD({\mathfrak a}):=e^X\CH^1$ equipped with the norm
\begin{equation*}
\|u\|_{\mathfrak a}^2:=\|u\|_{L^2}^2+\|e^{-X}u\|_{\CH^1}^2.
\end{equation*}
\end{definition}

\medskip

Since $e^X\in\CC^{1-\kappa}$ for any $\kappa>0$, the domain $\CD(\mathfrak{a})$ is dense in $\CH^{1-\kappa}$ thus in $L^2$. The following proposition states that this densely defined form is continuous and bounded from below.

\medskip

\begin{proposition}
There exists a random constant $C>0$ such that
\begin{equation*}
|{\mathfrak a}(u_1,u_2)|\le C\|u_1\|_{\mathfrak a} \|u_2\|_{\mathfrak a}
\end{equation*}
for $u_1,u_2\in\CD(\mathfrak{a})$. The form ${\mathfrak a}$ is quasi-coercive, i.e., there exists  random constants  $ \delta, C'>0$ such that
\begin{equation*}
{\mathfrak a}(u,u) + C' \|u\|_{L^2}^2 \ge \delta \|u\|_{\mathfrak a}^2 
\end{equation*}
for all $u=e^Xv\in\CD({\mathfrak a})$.
\end{proposition}

\medskip

\begin{proof}
The continuity follows directly from
\begin{equation*}
\big|\langle|\nabla X|^{2\diamond}e^{2X},v_1v_2\rangle\big|\le\||\nabla X|^{2\diamond}\|_{\CC^{-\kappa}}\|e^{2X}\|_{\CC^{1-\kappa}}\|v_1\|_{\CH^1}\|v_2\|_{\CH^1}.
\end{equation*}
Now we prove the second statement. Set $u = e^X v$ with $v \in \CH^1$. 
We have for any $\kappa >0$
\begin{align*}
&\hspace{-1cm}{\mathfrak a}(u,u) - \xi_0\int_{\IT^2}|v(x)|^2e^{2X(x)}\drm x \\
&=\int_{\IT^2} |\nabla v(x) |^2(x)e^{2X(x)}\drm x -\langle|\nabla X|^{2\diamond}v,ve^{2X}\rangle\\
&\ge e^{-\| X\|_{L^\infty}}  \int_{\IT^2} |\nabla v(x) |^2(x)\drm x -\||\nabla X|^{2\diamond}\|_{\CC^{-\kappa}}\|e^{2X}\|_{\CC^{2\kappa}}\|v\|_{\CH^{2\kappa}}^2. 
\end{align*}
For small $\kappa > 0$ we use the standard interpolation inequality, which is valid for every $\varepsilon > 0$,
\begin{equation*}
\|v\|_{\CH^{2\kappa}} \le \varepsilon \|v\|_{\CH^{1}} + c_\varepsilon \|v\|_{L^2}
\end{equation*}
for some $c_\varepsilon > 0$. We choose $\varepsilon$ small enough and insert this inequality in the previous estimates to obtain the statement.  
\end{proof}
As a consequence of the previous proposition one obtains that the norms $\|\cdot\|_{\CD(\mathfrak{a})}$ and $\|e^{-X}\cdot\|_{\CH^1}$ are equivalent. 

We now prove that the form is closed.

\medskip

\begin{proposition}
The form $\mathfrak{a}$ is closed, that is $(\CD(\mathfrak{a}),\|\cdot\|_\mathfrak{a})$ is a complete space.
\end{proposition}

\medskip

\begin{proof}
Let $(u_n)_{n\ge0}\subset\CD(\mathfrak{a})$ be a Cauchy sequence. Then $(e^{-X}u_n)_{n\ge0}$ is a Cauchy sequence in $\CH^1$ thus converges to a limit $v\in\CH^1$ while $(u_n)_{n\ge0}$ is a Cauchy sequence in $L^2$ thus converges to $u\in L^2$. We have
\begin{align*}
\|u-e^Xv\|_{L^2}&\le\|u-u_n\|_{L^2}+\|u_n-e^Xv\|_{L^2}\\
&\le\|u-u_n\|_{L^2}+\|e^X\|_{L^\infty}\|e^{-X}u_n-v\|_{L^2}
\end{align*}
hence $u=e^Xv\in\CD(\mathfrak{a})$ and this completes the proof. 
\end{proof}

Finally, we prove that $\mathfrak{a}$ is the limit in some sense of the renormalized form $\mathfrak{a}_\eps$.

\medskip

\begin{proposition}
For any $\kappa>0$, there exists a constant $C>0$ such that
\begin{equation*}
\big|\mathfrak{a}(u_1,u_2)-\mathfrak{a}_\eps(u_1^\eps,u_2^\eps)\big|\le C\|\Xi-\Xi_\eps\|_{\CX^\kappa(\IT^2)}\|v_1\|_{\CH^1}\|v_2\|_{\CH^1}
\end{equation*}
with $u_i^\eps=e^{X_\eps}v_i$ for $\eps\ge0$.
\end{proposition}

\medskip

\begin{proof}
Let $v_1,v_2\in\CH^1$ and consider $u_i^\eps=e^{X_\eps}v_i$ for $\eps\ge0$, i.e., $u_i^\eps\in\CH^1$ the form domain of $a_\eps$ for any $\eps>0$ while $u_i\in e^X\CH^1$ for $\eps=0$. We have
\begin{equation*}
\mathfrak{a}_\eps(u_1^\eps,u_2^\eps)=\langle \nabla v_1,\nabla v_2\rangle_{L^2(\IT^2,e^{2X}\drm x)}-\big\langle|\nabla X_\eps|^2-c_\eps,v_1v_2e^{2X_\eps}\big\rangle+\langle\xi_\eps,1\rangle\langle v_1,v_2\rangle_{L^2(\IT^2,e^{2X_\eps}\drm x)}
\end{equation*}
and hence 
\begin{align*}
\big|\mathfrak{a}(u_1,u_2)-\mathfrak{a}_\eps(u_1^\eps,u_2^\eps)\big|&\le\big|\langle|\nabla X_\eps|^2-c_\eps-|\nabla X|^{2\diamond},v_1v_2\rangle\big|+\big|\big\langle\langle\xi_\eps,1\rangle e^{X_\eps}-\xi_0e^X,v_1v_2\big\rangle\big|\\
&\le C\|\Xi-\Xi_\eps\|_{\CX^\kappa}\|v_1\|_{\CH^1}\|v_2\|_{\CH^1}
\end{align*}
for any $\kappa>0$ and the proof is complete.
\end{proof}

\section{Construction in three dimensions} \label{Section3D}

In three dimensions, the expression
\begin{equation*}
\langle \nabla v_1,\nabla v_2\rangle_{L^2(\IT^3,e^{2X}\drm x)}-\big\langle|\nabla X|^{2\diamond},v_1v_2e^{2X}\big\rangle+\xi_0\langle v_1,v_2e^{2X}\rangle
\end{equation*}
does not make sense anymore for $v_1,v_2\in\CH^1$ since $|\nabla X|^{2\diamond}$ belongs to $\CC^{-1-\kappa}$ for any $\kappa>0$. In this case, one makes the change of variable $u=e^{X+Y}v$ with $Y$ the solution to
\begin{equation*}
\Delta Y=|\nabla X|^{2\diamond}-\big\langle|\nabla X|^{2\diamond},1\big\rangle
\end{equation*}
which belongs to $\CC^{1-\kappa}$. We have
\begin{equation*}
Hu=-e^{X+Y}\Delta v-2e^{X+Y}(\nabla X+\nabla Y)\cdot\nabla v-\big(|\nabla Y|^{2\diamond}+2\nabla X\diamond\nabla Y-\langle |\nabla X|^{2\diamond},1\rangle-\xi_0\big)e^{X+Y}v
\end{equation*}
hence
\begin{align*}
\mathfrak{a}(u_1,u_2)&=\langle Hu_1,u_2\rangle\\
&=\langle He^{X+Y}v_1,e^{X+Y}v_2\rangle\\
&=-\langle \Delta v_1,v_2e^{2X+2Y}\rangle-2\langle \nabla(X+Y)\cdot\nabla v_1,v_2e^{2X+2Y}\rangle-\langle|\nabla Y|^{2\diamond},v_1v_2e^{2X+2Y}\rangle\\
&\quad-2\langle\nabla X\diamond\nabla Y,v_1v_2e^{2X+2Y}\rangle+\langle|\nabla X|^{2\diamond},1\rangle+\xi_0,v_1v_2e^{2X+2Y}\rangle\\
&=\int_{\IT^3}\nabla v_1(x)\cdot\nabla v_2(x)e^{2X(x)+2Y(x)}\drm x-\langle|\nabla Y|^{2\diamond}+2\nabla X\diamond\nabla Y,v_1v_2e^{2X+2Y}\rangle\\
&\quad+\big(\langle|\nabla X|^{2\diamond},1\rangle+\xi_0\big)\int_{\IT^3}v_1(x)v_2(x)e^{2X(x)+2Y(x)}\drm x
\end{align*}
which is well-defined for $v_1,v_2\in\CH^1$ since $|\nabla Y|^{2\diamond}\in\CC^{-\kappa}$ and $\nabla X\diamond\nabla Y\in\CC^{-\frac{1}{2}-\kappa}$ for any $\kappa>0$.

\medskip

\begin{definition}
The Anderson form is defined by 
\begin{align*}
\mathfrak{a}(u_1,u_2)&:=\langle \nabla v_1,\nabla v_2\rangle_{L^2(\IT^3,e^{2X+2Y}\drm x)}-\langle|\nabla Y|^{2\diamond}+2\nabla X\diamond\nabla Y,v_1v_2e^{2X+2Y}\rangle\\
&\quad+(\langle|\nabla X|^{2\diamond},1\rangle+\xi_0)\langle v_1,v_2\rangle_{L^2(\IT^3,e^{2X+2Y}\drm x)}
\end{align*}
where $v_i=e^{-X}u_i$ with domain $\CD(\mathfrak{a}):=e^{X+Y}\CH^1$ equipped with the norm
\begin{equation*}
\|u\|_\mathfrak{a}^2:=\|u\|_{L^2}^2+\|e^{-(X+Y)}u\|_{\CH^1}^2.
\end{equation*}
\end{definition}

Since $e^{X+Y}\in\CC^{\frac{1}{2}-\kappa}$ for any $\kappa>0$, the domain $\CD(\mathfrak{a})$ is dense $\CH^{\frac{1}{2}-\kappa}$ thus in $L^2$. The following proposition states that this densely defined form is continuous and bounded from below. The proofs are obtained following the same path as in two dimensions.

\medskip

\begin{proposition}
There exists a random constant $C>0$ such that
\begin{equation*}
|\mathfrak{a}(u_1,u_2)|\le C\|u_1\|_\mathfrak{a}\|u_2\|_\mathfrak{a}
\end{equation*}
for $u_1,u_2\in\CD(\mathfrak{a})$. There exists random constants $\delta, C'>0$ such that
\begin{equation*}
{\mathfrak a}(u,u) + C' \|u\|_{L^2}^2 \ge \delta \|u\|_{\mathfrak a}^2
\end{equation*}
for all $u=e^{X+Y}v\in\CD(\mathfrak{a})$.
\end{proposition}
Again, as in the two dimension case the form is closed.

\medskip

\begin{proposition}
The form $\mathfrak{a}$ is closed, that is $(\CD(\mathfrak{a}),\|\cdot\|_\mathfrak{a})$ is a complete space.
\end{proposition}

\medskip

Finally, $\mathfrak{a}$ is also the limit in some sense of the renormalized form $\mathfrak{a}_\eps$.

\medskip

\begin{proposition}
For any $\kappa>0$, there exists a constant $C>0$ such that
\begin{equation*}
\big| \mathfrak{a}(u_1,u_2)- \mathfrak{a}_\eps(u_1^\eps,u_2^\eps)\big|\le C\|\Xi-\Xi_\eps\|_{\CX^\kappa(\IT^3)}\|v_1\|_{\CH^1}\|v_2\|_{\CH^1}
\end{equation*}
with $u_i^\eps=e^{X_\eps+Y_\eps}v_i$ for $\eps\ge0$.
\end{proposition}

\section{Positivity and spectral gap} \label{SectionSpectralGap}

The construction of the form $\mathfrak{a}$  is the same in two and three dimensions. It  is densely defined, symmetric bounded from below, continuous and closed. Its associated operator $H$ has domain
\begin{equation*}
\CD(H)=\big\{u\in L^2(\IT^d)\ ;\ \exists v\in L^2(\IT^d),\forall\varphi\in\CD(\mathfrak{a}), \ \mathfrak{a}(u,\varphi)=\langle v,\varphi\rangle\big\}.
\end{equation*}
The operator $H$ is self-adjoint, densely defined and bounded from below. Since $\CD(\mathfrak{a})$ is imbedded  into a Sobolev space of positive regularity, it is compactly imbedded  in $L^2(\IT^d)$. Therefore, $H$  has discrete spectrum 
\begin{equation*}
\lambda_1\le\lambda_2\le\ldots
\end{equation*}
and there exists an orthonormal basis of $L^2(\IT^d)$ which is given by  eigenfunctions of $H$. An important information  is the existence of a spectral gap with a positive  ground state.  This is already known (see for example \cite{BailleulDangMouzard22}) and it is a key   to prove two-sided Gaussian bounds for the corresponding heat kernel of $H$. By the classical Krein-Rutman theorem, the  general idea to get a spectral gap with a positive ground state is to prove that the semigroup $e^{-tH}$ is positive and irreducible. This means that for any non-negative (and nontrivial) $f \in L^2(\IT^d)$, we have at any time $t > 0$, $e^{-tH} f > 0$ a.e. on $\IT^d$. The irreducibility is sometimes called {\it strict positivity} or {\it positivity improving}. 
Unlike \cite{BailleulDangMouzard22} which relies on quantitative study of the linear Parabolic Anderson Model equation and an approximation argument, we can obtain positivity and irreducibility readily from the form. These two properties are indeed characterized in terms of the form. See Theorems $2.6$ and  $2.10$ in \cite{Ouhabaz05}. Thus, we provide a pedestrian approach to the existence of a spectral gap even in three dimensions which usually relies on involved computations with expansion of order $5$ using regularity structures or paracontrolled calculus.

\medskip

\begin{theorem}
The semigroup $e^{-tH}$ is irreducible. In particular, the first eigenvalue is simple, that is $\lambda_1<\lambda_2$ and there exists  a  positive ground state $\Psi\in\CD(H)$.
\end{theorem}

\medskip

\begin{proof}
Both positivity and irreducibility are not changed under multiplication by $e^{X}$ or $e^{X+Y}$ and so we use the form 
$\mathfrak{a}$ constructed in the previous sections. 

Let $u \in D(\mathfrak{a})$ and $v \in \CH^1$ such that $u = e^{X} v$ if $d=2$ and $u = e^{X+Y} v$ if $d=3$. Then clearly, $u^+ = e^{X} v^+$ (or $e^{X+Y} v^+$) and $u^- = e^X v^-$ (or $e^{X+Y} v^-$). Since $v^+, v^- \in \CH^1$, we have $u^+, u^- \in D(\mathfrak{a})$. In addition, it is obviously seen from the definition of the Anderson form that $ \mathfrak{a}(u^+, u^-) = 0$. By \cite{Ouhabaz05}, Theorem $2.6$ we conclude that $(e^{-tH})_{t\ge 0}$ is a positive semigroup. 

Now we prove irreducibility. We apply Theorem $2.10$ from \cite{Ouhabaz05}. Since $H$ is a local operator, it is enough to prove that if
$D\subset\IT^d$ is such that
\begin{equation*}
\forall u\in D(\mathfrak{a}),\quad \IDC_D u\in D(\mathfrak{a}),
\end{equation*}
then either  $|D|=0$ or $|\IT^d\backslash D|=0$. 
Clearly, $ \IDC_D u = e^{X}  \IDC_D v $ if $d=2$ and $ \IDC_D u = e^{X +Y}  \IDC_D v $ if $d= 3$. This 
implies
\begin{equation*}
\forall v\in\CH^1,\quad\IDC_Dv\in\CH^1. 
\end{equation*}
Theorem $2.10$ from \cite{Ouhabaz05} applied to the Laplacian (whose form domain is $\CH^1$) gives that 
$|D|=0$ or $|\IT^d\backslash D|=0$. This proves irreducibility.

The rest of the theorem is classical and it is a direct consequence of the Krein-Rutman theorem. 
\end{proof}

\bibliographystyle{siam}
\bibliography{biblio.bib}

\vspace{2cm}

\noindent \textcolor{gray}{$\bullet$} CNRS \& Department of Mathematics and Applications, ENS Paris, 45 rue d'Ulm, 75005 Paris, France\\
{\it E-mail}: antoine.mouzard@math.cnrs.fr

\noindent \textcolor{gray}{$\bullet$} E. Ouhabaz --  Université de Bordeaux, CNRS, Institut de Mathématiques, F-33405 Bordeaux, France.\\
{\it E-mail}: elmaati.ouhabaz@math.u-bordeaux.fr

\end{document}